\documentclass{amsart}

\newtheorem{theorem}{Theorem}[section]

\theoremstyle{definition}

\newtheorem{example}[theorem]{Example}

\theoremstyle{remark}

\begin{document}

\title[Trigonometric integrals]{The integrals 
in Gradhteyn and Ryzhik. \\
 Part 5: Some 
trigonometric integrals}

\author{Tewodros Amdeberhan}
\address{Department of Mathematics,
Tulane University, New Orleans, LA 70118}
\email{tamdeberhan@math.tulane.edu}

\author{Luis A. Medina}
\address{Department of Mathematics,
Tulane University, New Orleans, LA 70118}
\email{lmedina@math.tulane.edu}

\author{Victor H. Moll}
\address{Department of Mathematics,
Tulane University, New Orleans, LA 70118}
\email{vhm@math.tulane.edu}

\subjclass{Primary 33}

\date{\today}


\begin{abstract}
We present evaluations and provide proofs of definite integrals involving the 
function $x^{p} \cos^{n} x$. These formulae are generalizations of 
$3.761.11$ and $3.822.1$, among others, in the classical table of integrals by 
I. S. Gradshteyn and I. M. Ryzhik.
\end{abstract}

\maketitle

\textwidth=6in
\textheight=8.5in

\newcommand{\nn}{\nonumber}
\newcommand{\ba}{\begin{eqnarray}}
\newcommand{\ea}{\end{eqnarray}}
\newcommand{\ift}{\int_{0}^{\infty}}
\newcommand{\no}{\noindent}
\newcommand{\at}{\text{tan}^{-1}}
\newcommand{\X}{{\mathbb{X}}}
\newcommand{\realpart}{\mathop{\rm Re}\nolimits}
\newcommand{\imagpart}{\mathop{\rm Im}\nolimits}

\newtheorem{Definition}{\bf Definition}[section]
\newtheorem{Thm}[Definition]{\bf Theorem} 
\newtheorem{Lem}[Definition]{\bf Lemma} 
\newtheorem{Cor}[Definition]{\bf Corollary} 
\newtheorem{Prop}[Definition]{\bf Proposition} 
\newtheorem{Note}[Definition]{\bf Note}
\numberwithin{equation}{section}

\section{Introduction} \label{intro} 
\setcounter{equation}{0}

The table of integrals \cite{gr} contains a large variety of 
evaluations of the type
\begin{equation}
I = \int_{a}^{b} A(x) R(\sin x, \cos x) \, dx
\end{equation}
\noindent
where $A$ is an algebraic function, $R$ is rational and 
$ - \infty \leq a < b \leq \infty$. We 
present a systematic discussion of two families of integrals of 
this type. This paper is part of a general program started 
in \cite{moll-gr1,moll-gr2,moll-gr3,moll-gr4} intended to 
provide proofs and context 
to the formulas in \cite{gr}. 

The first class considered here corresponds to the complete integrals
\begin{equation}
c(n,p) := \int_{0}^{\pi/2} x^{p} \cos^{n}x \, dx, 
\end{equation}
\noindent
and 
\begin{equation}
s(n,p) := \int_{0}^{\pi/2} x^{p} \sin^{n}x \, dx, 
\end{equation}
\noindent
where $n, \, p \in \mathbb{N}$. In section \ref{first} we present closed-form
expressions for these integrals. 
These expressions involve the sums 
\begin{equation}
\sum_{1 \leq k_{1} \leq k_{2} \leq \cdots \leq k_{j} \leq n} 
\frac{1}{k_{1}^{2} k_{2}^{2} \cdots k_{j}^{2}}, \label{eulersums}
\end{equation}
\noindent
that are closely related to the  multiple zeta values
\begin{equation}
\zeta(i_{1},i_{2}, \ldots, i_{s}) = \sum_{0 < k_{1} < k_{2} < \cdots< k_{s}} 
\frac{1}{k_{1}^{i_{1}} k_{2}^{i_{2}} \cdots k_{s}^{i_{s}} }.
\end{equation}
\noindent
The reader will find in Section 3.4 of \cite{borw2} an introduction to 
these sums. 

In general, one does not  expect
such elementary evaluations to extend to 
$p \not \in \mathbb{N}$. For example, the change 
of variables $x = \pi t^{2}/2$ produces 
\begin{equation}
\int_{0}^{\pi/2} x^{-1/2} \cos x \, dx = \sqrt{2 \pi} 
\int_{0}^{1} \cos \left( \frac{\pi t^{2}}{2} \right) \, dt. 
\end{equation}
\noindent
The latter integral is evaluated in terms of the {\em cosine Fresnel}
function
\begin{equation}
\text{FresnelC}[x]  := 
\int_{0}^{x} \cos \left( \frac{\pi t^{2}}{2} \right) \, dt,
\end{equation}
\noindent
which indeed is not an elementary function. 

The second class considered here presents generalizations of 
the formula $3.822.1$ in \cite{gr} stated as 
\begin{equation}
\ift \frac{\cos^{2n+1}x}{\sqrt{x}} \, dx = \frac{1}{2^{2n}} 
\sqrt{\frac{\pi}{2}} \sum_{k=0}^{n} \binom{2n+1}{n+k+1} \frac{1}{\sqrt{2k+1}}, 
\quad n \in \mathbb{N}.
\label{ex1}
\end{equation}
\noindent
The integral in (\ref{ex1}) can be transformed via $t = x^{2}$ 
to provide the evaluation of 
\begin{equation}
\ift \cos^{2n+1}t^{2} \, dt,
\end{equation}
\noindent
that is given as the case $p=2$ in Theorem \ref{thmluis1}. 

Section \ref{half-line} contains analytic expressions for the 
generalizations
\begin{equation}
C_{n}(p,b) := \ift x^{-p} \cos^{2n+1}(x+b) \, dx, 
\end{equation}
\noindent
and
\begin{equation}
S_{n}(p,b) := \ift x^{-p} \sin^{2n+1}(x+b) \, dx.
\end{equation}
\noindent
The last section also contains some evaluations obtained by differentiation 
with respect to parameters. An illustrative example is 
\begin{equation}
\ift \ift \frac{\log x \, \log y }{\sqrt{xy}} \cos(x+y) \, dx \, dy = 
( \gamma + 2 \log 2) \pi^{2},
\end{equation}
\noindent
that is equivalent to 
\begin{equation}
\ift \ift \log x \, \log y \,  \cos(x^{2}+y^{2}) \, dx \, dy = 
\frac{1}{16}( \gamma + 2 \log 2) \pi^{2}. 
\end{equation}
\noindent
A generalization of this evaluation appears as Example \ref{last-example}.

The method described in the present work gives impetus to a class of
integrals that are closely related to the particular integral 
computations addressed in this paper. 

\section{The first example} \label{first} 
\setcounter{equation}{0}

In this section we present the evaluation in closed-form of the 
definite integrals 
\begin{equation}
c(n,p) := \int_{0}^{\pi/2} x^{p} \cos^{n}x \, dx. 
\end{equation}
\noindent
A special case of this appears as $3.822.1$ in \cite{gr}.

The first step towards the evaluation of $c(n,p)$ is to produce a recurrence. 

\begin{Thm}
The integral $c(n,p)$ satisfies the recurrence 
\ba
c(n,p) & = & \frac{n-1}{n} c(n-2,p) - \frac{p(p-1)}{n^{2}} c(n,p-2), 
\label{rec1}
\ea
\noindent
for  $n \geq 2, \, p \geq 2$.
\end{Thm}
\begin{proof}
The identity $\cos^{2}x = 1 - \sin^{2}x$ yields 
\ba
c(n,p) & = & c(n-2,p) - I(n,p) \label{red1}
\ea
\no
where 
\ba
I(n,p) & := & \int_{0}^{\pi/2} x^{p} \cos^{n-2}x \,  \sin^{2} x \, dx. \nn
\ea
\no
Now 
\ba
I(n,p) & = & \int_{0}^{\pi/2} x^{p} \sin x \times \frac{d}{dx} \left( 
- \frac{1}{n-1} \cos^{n-1} x \right) \, dx \nn \\
   & = & \frac{1}{2n-1} \int_{0}^{\pi/2}  \left( x^{p} \cos x + 
                       p x^{p-1} \sin x \right) \, \cos^{n-1}x \, dx \nn \\
   & = & \frac{c(n,p)}{n-1} + \frac{p}{n-1} \int_{0}^{\pi/2} x^{p-1} \sin x 
\cos^{n-1} x \; dx. \nn
\ea
\no
Moreover
\ba
\int_{0}^{\pi/2} x^{p-1} \sin x \cos^{n-1} x \; dx & = & 
\int_{0}^{\pi/2} x^{p-1} 
\frac{d}{dx} \left( - \frac{1}{n} \cos^{n}x 
\right) \; dx \nn \\
  & = & \frac{p-1}{n} c(n,p-2). \nn
\ea
\end{proof}

\noindent
{\bf Strategy}: According to (\ref{rec1}), the 
integral $c(n,p)$ can be evaluated in terms of 
the initial values given in the table. The indices $m$ and $q$ have the 
same parity as $n$ and $p$ respectively and range over $0 \leq m \leq n$ and 
$0 \leq q \leq p$.  \\

\begin{center}
\begin{tabular}{||c||c||c||}
\hline
$n$ \text{ modulo } $2$  & $p$ \text{ modulo}  $2$& \text{initial conditions} \\
\hline 
$0$ & $0$ & $c(m,0)$ \, $c(0,q)$ \\
$1$ &  $0$  & $c(m,0)$ \, $c(1,q)$  \\
$0$ &  $1$  & $c(m,1)$ \, $c(0,q)$  \\
$1$ &  $1$  & $c(m,1)$ \, $c(1,q)$   \\
\hline
\end{tabular}
\end{center}

\medskip 

We now evaluate the initial conditions $c(n,0), \, c(n,1), \, c(0,p)$ and 
$c(1,p)$. \\

\noindent
{\bf The expression for} $c(0,p)$. 

The computation of the identity 
\begin{equation}
c(0,p) = \frac{1}{p+1} \left( \frac{\pi}{2} \right)^{p+1}
\label{c0p}
\end{equation}
\noindent
is immediate.


\medskip

\noindent
{\bf The expression for} $c(n,0)$.

\noindent
This is classical. The 
result appears as 
$3.621.3$ and $3.621.4$ in \cite{gr}.

\begin{Thm}
\label{thm-wallis}
{\bf (Wallis' formula and companion)}. Let $n \in {\mathbb{N}}_{0}$. Then 
\begin{equation}
\label{one1}
c(2n,0)  = 
 \frac{\pi}{2^{2n+1}} \binom{2n}{n},
\end{equation}
\noindent
and 
\begin{equation}
\label{one2}
c(2n+1,0) =  \frac{2^{2n}}{(2n+1) \binom{2n}{n}}.
\end{equation}
\end{Thm}

\medskip

The shortest proof of Theorem \ref{thm-wallis} employs the representation 
\begin{equation}
\label{one3}
c(n,0)  = \int_{0}^{\pi/2} \cos^{n}x \, dx = 
2^{n-1} B \left( \frac{n+1}{2}, \frac{n+1}{2} \right), 
\end{equation}
\noindent 
that appears as $3.621.1$ in \cite{gr}. Here $B$ is the {\em Euler's beta 
function} defined by the integral
\begin{equation}
B(x,y) = \int_{0}^{1} t^{x-1} (1-t)^{y-1} \, dt.
\end{equation}
\noindent
The expression (\ref{one3}) follows from the change of variables 
$t = \cos u$. To express (\ref{one1}) and (\ref{one2}), in terms of the 
beta function,  employ the standard relation 
\begin{equation}
B(x,y) = \frac{\Gamma(x) \, \Gamma(y)}{\Gamma(x+y)}, 
\end{equation}
\noindent
and the special values 
\begin{equation}
\Gamma(n) = (n-1)! \, \text{ and } 
\, \Gamma( n + \tfrac{1}{2} ) =  \frac{\sqrt{\pi} \, (2n)!}{2^{2n} \, n!} 
\end{equation}
\noindent
that are valid for $n \in \mathbb{N}$. \\

The identity in Theorem \ref{thm-wallis}, in the case $n$ is even, that is,
\begin{equation}
c(2n,0) = \int_{0}^{\pi/2} \cos^{2n} \theta \, d \theta 
 =  \frac{\pi}{2^{2n+1}} \binom{2n}{n}, \label{wallis12}
\end{equation}
\no
is Wallis's formula and sometimes found in calculus 
books (see e.g. \cite{lhe}, page 492). To prove it, first  
write $\cos^{2} \theta = 1 - \sin^{2} \theta$
and use integration by parts to obtain the recursion
\ba
c(2n,0) & = & \frac{2n-1}{2n} c(2n-2,0). \label{recur1}
\ea
\no
Then verify that the right side of (\ref{wallis12}) satisfies the 
same recurrence together with the initial value $\pi/2$ for $n=0$. 

We now present a new proof of Wallis's formula (\ref{wallis12})
in the context of rational integrals. Extensions of the ideas in
this proof have produced {\em rational Landen transformations}. The 
reader will find in \cite{boros2, boros1,hubbard1,manna-moll2,manna-moll1}
details on these transformations. 

Start with
\ba
c(2n,0)  & = & \int_{0}^{\pi/2} \cos^{2n} \theta \, d \theta  = 
      \int_{0}^{\pi/2} \left( \frac{1  + \cos 2 \theta}{2} \right)^{n} 
\, d \theta.  \nn  
\ea
\no
Now introduce $\psi = 2 \theta$ and expand and simplify the result by 
observing that the odd powers of cosine integrate to zero. The inductive 
proof of (\ref{wallis12}) requires 
\ba
c(2n,0) & = & 2^{-n} \sum_{i=0}^{\lfloor{n/2\rfloor}} \binom{n}{2i} c(2i,0). \label{recur}
\ea
\no
Note that $c(2n,0)$ is uniquely determined by (\ref{recur}) along with the
initial value $c(0,0) = \pi/2$. Thus  (\ref{wallis12}) now follows from the
identity 
\ba
f(n) := \sum_{i=0}^{\lfloor{n/2\rfloor}} 2^{-2i} \binom{n}{2i} \binom{2i}{i} & = &
2^{-n} \binom{2n}{n}. \label{sum1}
\ea

We now provide a mechanical proof of (\ref{sum1}) using the theory  
developed by Wilf and Zeilberger, which is explained in
\cite{nemes,aequalsb}; the sum in (\ref{sum1}) is the example used 
in \cite{aequalsb} (page 113) to 
illustrate their method. The command
$$ct(\text{binomial}(n,2i) \, \text{binomial}(2i,i) 2^{-2i}, 1, i, n,N)$$
produces
\ba
f(n+1) & = & \frac{2n+1}{n+1} \; f(n), \label{recur2}
\ea
\no
and one checks that $2^{-n} \binom{2n}{n}$ satisfies this 
recursion. Note that (\ref{recur1}) and (\ref{recur2}) are 
equivalent under
\ba
c(2n,0) & = & \frac{\pi}{2^{n+1}} f(n). \nn
\ea
\noindent
The proof is complete.

\medskip

\noindent
{\bf Closed form expression for} $c(1,p)$. 

We now consider the evaluation of 
\begin{equation}
c(1,p) := \int_{0}^{\pi/2} x^{p} \cos x \, dx.
\label{conep}
\end{equation}
\noindent
The following evaluation appears as $3.761.11$ in \cite{gr}.

\begin{Thm}
\label{thm24}
Let $p \in \mathbb{N}$ and $\delta_{\text{odd},p}$ be Kronecker's delta 
function at the odd integers. Then
\begin{equation}
c(1,p) = \sum_{k=0}^{\xi_{p}} (-1)^{k} 
\frac{p!}{(p-2k)!} \left( \frac{\pi}{2} \right)^{p-2k} 
- (-1)^{\xi_{p}} \delta_{\text{odd},p} \, p! 
\label{c1p}
\end{equation}
\noindent
where $\xi_{p} = \lfloor{ \frac{p}{2} \rfloor}$. 
\end{Thm}
\begin{proof}
Both sides of the equation (\ref{c1p}) satisfy the initial value problem
\begin{equation}
u_{p} - p(p-1)u_{p-2} = \left( \frac{\pi}{2} \right)^{p} \text{ and } \,
u_{0} = 1,  \, u_{1} = \frac{\pi-2}{2}. 
\label{recu00}
\end{equation}
\noindent
Actually the recurrence (\ref{recu00}) is obtained  using integration
by parts in (\ref{conep}). Iterating this recurrence yields the 
right hand side of (\ref{c1p}). 
\end{proof}

\begin{Note}
The result in Theorem \ref{thm24} can be expressed in terms of 
the Taylor polynomial for $\cos x$:
\begin{equation}
f_{p}(x) = (-1)^{\xi_{p}} p! \left( -1 + 
\sum_{k=0}^{\xi_{p} } \frac{(-1)^{k}}{(2k+1)! } x^{2k+1} 
\right).
\end{equation}
\noindent
The formula (\ref{c1p}) can be restated 
\begin{equation}
c(1,p) = \begin{cases}
          f_{p}(\pi/2), \text{ for } p \text{ odd}, \\ 
          f_{p}'(\pi/2), \text{ for } p \text{ even}. 
        \end{cases}
\end{equation}
\end{Note}

\medskip

\noindent
{\bf Closed form expression for} $c(n,1)$: in fact, this would be the 
last initial condition we require to execute the startegy outlined at the 
beginning of this section. 

\begin{Thm}
The integral $c(n,1)$ satisfies the recurrence
\ba
c(n,1) & = & \frac{n-1}{n} c(n-2,1) - \frac{1}{n^{2}}. \label{rec2}
\ea
\end{Thm}
\begin{proof}
The identity $\cos^{2}x = 1 - \sin^{2}x$ yields 
\begin{equation}
c(n,1) = c(n-2,1) - J,
\end{equation}
\noindent
where 
\begin{equation}
J = \int_{0}^{\pi/2} x \sin^{2}x \, \cos^{n-2}x \, dx. 
\end{equation}
\noindent
Integration by parts leads to 
\begin{equation}
J = \frac{1}{n-1} \int_{0}^{\pi/2} ( \sin x + x \cos x ) \cos^{n-1}x \, dx.
\end{equation}
\noindent
This produces (\ref{rec2}).
\end{proof}

\noindent
The solution of (\ref{rec2}) yields  a closed-form formula for $c(n,1)$.

\begin{Thm}
\label{ctwon}
The integral $c(n,1)$  is given according to the parity of $n$, by 
\begin{equation}
c(2n,1)  = \frac{\binom{2n}{n}}{2^{2n+2}} \left( \frac{\pi^{2}}{2} - 
\sum_{k=1}^{n} \frac{2^{2k}}{k^{2} \, \binom{2k}{k}} \right),
\label{value2}
\end{equation}
\noindent
for even indices. For odd indices, we have 
\begin{equation}
c(2n+1,1) = \frac{2^{2n}}{(2n+1) \binom{2n}{n} } 
\left( \frac{\pi}{2} - \sum_{k=0}^{n} \frac{\binom{2k}{k}}{2^{2k} \, (2k+1)} \right). 
\label{value1}
\end{equation}
\end{Thm}

To establish this result we solve a more general recurrence than (\ref{rec2}). 

\begin{Lem}
\label{lemrec}
Let $a_{n}, \, b_{n}$ and $r_{n}$ be sequences with $a_{n}, \, 
b_{n} \neq 0$. Assume that
$z_{n}$ satisfies 
\begin{equation}
a_{n}z_{n} = b_{n} z_{n-1} + r_{n}, \, n \geq 1
\label{recu0}
\end{equation}
\noindent
with initial condition $z_{0}$. Then 
\begin{equation}
\label{valuex}
z_{n} = \frac{b_{1}b_{2} \cdots b_{n}}{a_{1}a_{2} \cdots a_{n}} 
\left( z_{0} + \sum_{k=1}^{n} \frac{a_{1}a_{2} \cdots a_{k-1}}
{b_{1} b_{2} \cdots b_{k}} r_{k} \right). 
\end{equation}
\end{Lem}
\begin{proof}
Introduce the integrating factor $d_{n}$ with the property that 
$d_{n}b_{n} = d_{n-1}a_{n-1}$. The recurrence (\ref{recu0}) becomes 
\begin{equation}
D_{n} - D_{n-1}  = d_{n}r_{n}, 
\end{equation}
\noindent
where $D_{n} = d_{n}a_{n}z_{n}$. Therefore, by telescoping,
\begin{equation}
D_{n} = D_{0} + \sum_{k=1}^{n} d_{k}r_{k},
\label{valueD}
\end{equation}
\noindent
with $D_{0} = d_{0}a_{0}z_{0}$. To find the integrating factor, observe that
\begin{equation}
\frac{d_{n}}{d_{n-1}} = \frac{a_{n-1}}{b_{n}}. 
\end{equation}
\noindent
Thus 
\begin{equation}
d_{n} = \frac{a_{0} a_{1} \cdots a_{n-1}}{b_{1} b_{2} \cdots b_{n}} d_{0}. 
\end{equation}
\noindent
Replacing in (\ref{valueD}) yields (\ref{valuex}). 
\end{proof}

\begin{Cor}
\label{cor-sys}
Let $n \in \mathbb{N}$ and assume that $z_{n}$ satisfies 
\begin{equation}
2nz_{n} = (2n-1)z_{n-1} + r_{n}, \, n \geq 1,
\end{equation}
\noindent
with the initial condition $z_{0}$. Let
$\lambda_{j} = 2^{2j}\binom{2j}{j}^{-1}$, then 
\begin{equation}
z_{n} = \frac{1}{\lambda_{n}} \left( z_{0} + \sum_{k=1}^{n} 
\frac{\lambda_{k} r_{k}}{2k} \right). 
\end{equation}
\end{Cor}
\begin{proof}
Use $a_{n} = 2n$ and $b_{n} = 2n-1$ in Lemma \ref{lemrec}.
\end{proof}

\medskip

We now apply Lemma \ref{lemrec}  on 
the recurrence (\ref{rec2}), repeated here for 
convenience to the reader,
\begin{equation}
c(n,1)  =  \frac{n-1}{n} c(n-2,1) - \frac{1}{n^{2}}. \nonumber
\end{equation}
\noindent
Observe that this recurrence splits naturally into even and odd 
branches. The value of $c(2n,1)$ is determined completely by $c(0,1)$, and 
$c(2n+1,1)$ by $c(1,1)$. Hence, there is no computational 
interaction between $c(2n,1)$ and 
$c(2n+1,1)$. Let
$x_{n} = c(2n,1)$ so that $x_{n}$ satisfies 
\begin{equation}
2n x_{n} = (2n-1) x_{n-1}  - \frac{1}{4n},
\end{equation}
\noindent
with the initial condition 
\begin{equation}
x_{0} = c(0,1) = \frac{\pi^{2}}{8}. 
\end{equation}
\noindent
Similarly, $y_{n} = c(2n+1,1)$, the odd component of $c(n,1)$,  satisfies 
\begin{equation}
(2n+1) y_{n} = 2n y_{n-1}  - \frac{1}{2n+1}
\end{equation}
\noindent
and the initial condition 
\begin{equation}
y_{0} = c(1,1) = \frac{\pi}{2}-1.  
\end{equation}

The expressions for $z_{n}$ in Lemma \ref{lemrec} yield the formulas for 
$c(2n,1)$ and also $c(2n+1,1)$ in Theorem \ref{ctwon}. The proof is complete. 

\medskip

\begin{Note}
The finite sums in (\ref{value2}) and (\ref{value1}) 
do not have closed-form, but it is a
classical result that, in the limit, 
\begin{equation}
\sum_{k=1}^{\infty} \frac{2^{2k}}{k^{2} \, \binom{2k}{k} } = 
\frac{\pi^{2}}{2} 
\end{equation}
\noindent
and 
\begin{equation}
\sum_{k=0}^{\infty} \frac{\binom{2k}{k} }{2^{2k} (2k+1)} = 
\frac{\pi}{2}. 
\end{equation}
\end{Note}

\noindent
\begin{Note}
Formula $3.821.3$ in \cite{gr} gives formulas equivalent to 
(\ref{value2}) and (\ref{value1}), respectively. 
\end{Note}

\medskip

Finally, we conclude 
this section by presenting ithe sought for closed form expression for the 
integral $c(n,p)$, for arbitrary $n, \, p \in \mathbb{N}$. The 
recurrence (\ref{rec1}), in the case of even 
indices $n$, becomes
\begin{equation}
2nX_{n}(p) = (2n-1)X_{n-1}(p) - \frac{p(p-1)}{2n} X_{n}(p-2)
\label{rec-even}
\end{equation}
\noindent
where $X_{n}(p) = c(2n,p)$. The initial value 
\begin{equation}
X_{0}(p) = \frac{1}{(p+1) 2^{p+1}} \pi^{p+1} 
\label{xzero}
\end{equation}
\noindent
given in (\ref{c0p}) and the recurrence (\ref{rec-even}) show the 
existence of rational numbers $a_{n,p,p+1-2j}$ such that 
\begin{equation}
X_{n}(p) = \sum_{j=0}^{\xi_{p}} a_{n,p,p+1-2j} \pi^{p+1-2j}, 
\end{equation}
\noindent 
with $\xi_{p} = \lfloor{ \frac{p}{2} \rfloor}$. The 
recurrence (\ref{rec-even}) is now expanded as 
\begin{eqnarray}
2n \sum_{j=0}^{\xi_{p}} a_{n,p,p+1-2j} \pi^{p+1-2j} & = & 
(2n-1) \sum_{j=0}^{\xi_{p}} a_{n-1,p,p+1-2j} \pi^{p+1-2j} \label{rec-all} \\
& - & \frac{p(p-1)}{2n} \sum_{j=0}^{\xi_{p-1}} a_{n,p-2,p-1-2j} \pi^{p-1-2j}.
\nonumber 
\end{eqnarray}

The fact that the coefficients $a_{n,p,j} \in \mathbb{Q}$ allows us to match 
the corresponding powers of $\pi$ in (\ref{rec-all}). The  highest order term
is $\pi^{p+1}$. Only two of the sums contain this power, therefore
\begin{equation}
2na_{n,p,p+1} = (2n-1) a_{n-1,p,p+1}. 
\label{recu-1}
\end{equation}
\noindent
The initial condition 
\begin{equation}
a_{0,p,p+1} = \frac{1}{(p+1)2^{p+1}}
\label{init-1}
\end{equation}
\noindent
comes from (\ref{xzero}). The solution to the initial value problem 
(\ref{recu-1}, \ref{init-1}) is then  found using 
Corollary \ref{cor-sys} (here $r_{n}=0$), namely that
\begin{equation}
a_{n,p,p+1} = \frac{\binom{2n}{n} }{(p+1)2^{2n+p+1}}. 
\label{formula1}
\end{equation}

The coefficient of the next highest power $\pi^{p-1}$, in 
(\ref{rec-all}), yields the recurrence 
\begin{equation}
2n a_{n,p,p-1} = (2n-1)a_{n-1,p,p-1} - \frac{p(p-1)}{2n} a_{n,p-2,p-1}. 
\label{recu-2}
\end{equation}
\noindent
Observe that the last term in this relation is given by (\ref{formula1}). 
Moreover, (\ref{xzero}) shows that $a_{0,p,p-1} = 0$. The solution to 
(\ref{recu-2}), following Corollary \ref{cor-sys}, is 
\begin{equation}
a_{n,p,p-1} = - \frac{p \binom{2n}{n} }{2^{2n+p+1}} \sum_{k_{1}=1}^{n} \frac{1}{k_{1}^{2}}. 
\label{formula2}
\end{equation}
\noindent
The next power of $\pi$ in (\ref{rec-all}) produces
\begin{equation}
2na_{n,p,p-3} = (2n-1)a_{n-1,p,p-3} + 
\frac{p(p-1)(p-2)}{n2^{2n+p}} \binom{2n}{n} \sum_{k_{1}=1}^{n} \frac{1}{k_{1}^{2}},
\end{equation}
\noindent
with $a_{0,p,p-3}=0$. One more use of Corollary \ref{cor-sys} yields 
\begin{equation}
a_{n,p,p-3} = \frac{\binom{2n}{n} \, p!}{2^{2n+p+1} \, (p-3)!} 
\sum_{k_{2}=1}^{n} \sum_{k_{1}=1}^{k_{2}} \frac{1}{k_{1}^{2}k_{2}^{2}}. 
\end{equation}

This procedure can be repeated until all descending powers of $\pi$ have been 
exhausted, hence a complete closed form for the integrals $c(n,p)$ will 
be made possible. 

\medskip

\begin{Thm}
Let $n, \, p \in \mathbb{N}$ and let 
$\xi_{p} = \lfloor{ \frac{p}{2} \rfloor}$.  Then the 
even branches $X_{n}(p) = c(2n,p)$ 
of the integral 
\begin{equation}
c(n,p) = \int_{0}^{\pi/2} x^{p} \cos^{n}x \, dx 
\end{equation}
\noindent
are given by
\begin{equation}
X_{n}(p) = \sum_{j=0}^{\xi_{p}}  a_{n,p,p+1-2j} \pi^{p+1-2j} 
+ \delta_{\text{odd},p} \cdot a_{n,p}^{*}, 
\end{equation}
\noindent
and the value of $a_{n,p,p+1-2j}$ 
for $p \geq 2$ 
and $0 \leq j \leq \xi_{p}$ 
is given by
\begin{equation}
a_{n,p,p+1-2j} =\frac{(-1)^{j} \binom{2n}{n} p!}{2^{2n+p+1} \, (p+1-2j)!} 
\sum_{1 \leq k_{1} \leq k_{2} \leq \cdots \leq k_{j} \leq n} 
\frac{1}{k_{1}^{2} k_{2}^{2} \cdots k_{j}^{2}}, \nonumber
\end{equation}
\noindent
and 
\begin{equation}
a_{n,p}^{*} =\frac{(-1)^{\xi_{p}} \binom{2n}{n} p!}{2^{2n}} 
\sum_{1 \leq k_{1} \leq k_{2} \leq \cdots \leq k_{p} \leq n} 
\frac{1}{k_{1}^{2} k_{2}^{2} \cdots k_{p}^{2}} 
\,  \sum_{j=1}^{k_{p}} \frac{2^{2j}}{j^{2} \binom{2j}{j}}. 
\end{equation}

Similarly, for the odd branches $Y(n,p) = c(2n+1,p)$ we have 
\begin{equation}
Y_{n}(p) = \sum_{j=0}^{\xi_{p}}  b_{n,p,p-2j} \pi^{p-2j} + 
\delta_{\text{odd},p} \cdot b_{n,p}^{*}, 
\end{equation}
\noindent
with 
\begin{equation}
b_{n,p,p-2j} =
\frac{(-1)^{j} \, p! \, 2^{2n+2j-p} }{(2n+1) \binom{2n}{n} \, (p-2j)!} 
\sum_{0 \leq k_{1} \leq k_{2} \leq \cdots \leq k_{j} \leq n} 
\frac{1}{(2k_{1}+1)^{2} (2k_{2}+1)^{2} \cdots (2k_{j}+1)^{2}}, \nonumber
\end{equation}
\noindent
and 
\begin{equation}
b_{n,p}^{*} = 
\frac{(-1)^{\xi_{p}} \, p! \, 2^{2n} }{(2n+1) \binom{2n}{n}} 
\sum_{0 \leq k_{1} \leq k_{2} \leq \cdots \leq k_{p} \leq n} 
\frac{1}{(2k_{1}+1)^{2} (2k_{2}+1)^{2} \cdots (2k_{p}+1)^{2}}
\, \sum_{j=0}^{k_{p}} \frac{\binom{2j}{j}}{2^{2j} (2j+1)}. \nonumber
\end{equation}
\end{Thm}

\medskip

\section{Some examples on the halfline} \label{half-line} 
\setcounter{equation}{0}

In this section we provide an analytic expression for 
\begin{equation}
C_{n}(p,b) = \ift x^{-p} \cos^{2n+1}(x+b) \, dx, 
\end{equation}
\noindent
and 
\begin{equation}
S_{n}(p,b) = \ift x^{-p} \sin^{2n+1}(x+b) \, dx.
\end{equation}

In the table \cite{gr} the evaluation of the special case $p = \frac{1}{2}$ 
and $b=0$:
\begin{equation}
\ift \frac{\cos^{2n+1}x}{\sqrt{x}} \, dx = \frac{1}{2^{2n}} 
\sqrt{\frac{\pi}{2}} \sum_{k=0}^{n} \binom{2n+1}{n+k+1} \frac{1}{\sqrt{2k+1}}, 
\label{lastone}
\end{equation}
\noindent
and 
\begin{equation}
\ift \frac{\sin^{2n+1}x}{\sqrt{x}} \, dx = \frac{1}{2^{2n}} 
\sqrt{\frac{\pi}{2}} \sum_{k=0}^{n} (-1)^{k} \binom{2n+1}{n+k+1} \frac{1}{\sqrt{2k+1}}, 
\label{lastone1}
\end{equation}
\noindent
as $3.822.2$ and $3.821.14$. 

\begin{Thm}
\label{main1}
Let $0 < p < 1$ and $n \in {\mathbb{N}}_{0} := \mathbb{N} \cup \{ 0 \}$. 
Then 
\begin{equation}
\ift x^{-p} \cos^{2n+1}x \, dx = \frac{\Gamma(1-p)}{2^{2n}} 
\sin \left( \frac{\pi p }{2} \right) 
\sum_{k=0}^{n} \frac{\binom{2n+1}{n-k}}{(2k+1)^{1-p}},
\end{equation}
\noindent
and 
\begin{equation}
\ift x^{-p} \sin^{2n+1}x \, dx = \frac{\Gamma(1-p)}{2^{2n}} 
\cos \left( \frac{\pi p }{2} \right) 
\sum_{k=0}^{n} (-1)^{k} \frac{\binom{2n+1}{n-k}}{(2k+1)^{1-p}}.
\end{equation}
\end{Thm}
\begin{proof}
The identity $2 \cos x = e^{ix} + e^{-ix}$ and the binomial theorem yield
\begin{equation}
\ift x^{-p} \cos^{2n+1} x \, dx = 
2^{-2n-1} \sum_{k=0}^{n} \binom{2n+1}{k} \ift 
x^{-p} \left( e^{i(2n+1-2k)x} + e^{-i(2n+1-2k)} \right) \, dx. 
\label{sum00}
\end{equation}
\noindent
Recall the Heaviside step function defined by $H(x) = 1$, if $x > 0$ and 
$H(x) = 0 $ otherwise. Then, each of the integrals in (\ref{sum00}) is 
evaluated using the Fourier transform
\begin{equation}
\int_{-\infty}^{\infty} H(x) x^{-p} e^{-i \omega x} \, dx = 
\frac{\Gamma(1-p)}{|\omega|^{1-p}} \text{exp}( - i \pi(1-p) 
\text{sign}(\omega)/2). 
\end{equation}
\end{proof}

\begin{Cor}
\label{thmluis1}
Let $p> 1$ be real and $n \in {\mathbb{N}}_{0}$. Then 
\begin{equation}
\ift \cos^{2n+1}x^{p} \, dx = 
\frac{1}{2^{2n}} \Gamma \left( \frac{p+1}{p} \right) \cos \left( \frac{\pi}{2p} \right)
 \sum_{k=0}^{n} \frac{\binom{2n+1}{n-k} }
{(2k+1)^{1/p}},
\label{cosq}
\end{equation}
\noindent
and
\begin{equation}
\ift \sin^{2n+1}x^{p} \, dx = 
\frac{1}{2^{2n}} \Gamma \left( \frac{p+1}{p} \right) \sin \left( \frac{\pi}{2p} \right)
 \sum_{k=0}^{n} (-1)^{k} \frac{\binom{2n+1}{n-k} }
{(2k+1)^{1/p}}.
\label{sineq}
\end{equation}
\end{Cor}
\begin{proof}
The change of variables $x \mapsto x^{1/(1-p)}$ in the results of Theorem 
\ref{main1} gives the result. 
\end{proof}

The last result described here is a further generalization of Theorem 
\ref{main1}. 

\begin{Thm}
\label{thmCS}
Assume $b \in \mathbb{R}, \, 0 < p < 1$ and $n \in {\mathbb{N}}_{0}$. Define
\begin{equation}
C_{n}(p,b) = \ift x^{-p} \cos^{2n+1}(x+b) \, dx 
\end{equation}
\noindent
and 
\begin{equation}
S_{n}(p,b) = \ift x^{-p} \sin^{2n+1}(x+b) \, dx.
\end{equation}
\noindent
Then 
\begin{equation}
C_{n}(p,b) = \frac{ \Gamma(1-p)}{2^{2n}}
\sum_{k=0}^{n} \binom{2n+1}{n-k} 
\frac{\sin( \tfrac{\pi p}{2} - (2k+1)b)}{(2k+1)^{1-p}}, 
\label{gencos}
\end{equation}
\noindent
and
\begin{equation}
S_{n}(p,b) = \frac{ \Gamma(1-p)}{2^{2n}}
\sum_{k=0}^{n} (-1)^{k} \binom{2n+1}{n-k} 
\frac{\cos( \tfrac{\pi p}{2} - (2k+1)b)}{(2k+1)^{1-p}}. 
\label{gensin}
\end{equation}
\end{Thm}
\begin{proof}
Denote the left-hand side of (\ref{gencos}) and (\ref{gensin}) by 
$f_{n}(b)$ and $g_{n}(b)$ respectively. Differentiation with respect to the
parameter $b$ yields
\begin{eqnarray}
\frac{\partial g_{n}}{\partial b} - (-1)^{n} (2n+1) f_{n} & = & 
(2n+1) \sum_{j=0}^{n-1} (-1)^{j} \binom{n}{j} f_{j}(b) \label{eone} \\
\frac{\partial f_{n}}{\partial b} + (-1)^{n} (2n+1) g_{n} & = & 
-(2n+1) \sum_{j=0}^{n-1} (-1)^{j} \binom{n}{j} g_{j}(b). \nonumber
\end{eqnarray}
\noindent
Considering $b$ and $p$ fixed, we now show that the right-hand side of 
(\ref{gencos}) and (\ref{gensin}) satisfy the system (\ref{eone}) with the
same initial conditions. This will establish the result. 

In the case of the right-hand side of (\ref{gencos}), it is required to 
check the identity
\begin{equation}
2^{-2n} \sum_{k=0}^{n} (-1)^{k} \binom{2n+1}{n-k} 
\frac{\sin( \tfrac{\pi}{2}p -(2k+1)b)}{(2k+1)^{1-p}} = \nonumber
\end{equation}
\begin{equation}
(2n+1) \sum_{j=0}^{n} (-1)^{j} \binom{n}{j} 2^{-2j} 
\sum_{k=0}^{j} \binom{2j+1}{j-k} 
\frac{\sin( \tfrac{\pi}{2}p -(2k+1)b)}{(2j+1)^{1-p}}.  \nonumber
\end{equation}
\noindent
To verify this we compare the coefficients of the transcendental terms
\begin{equation}
\frac{\sin( \tfrac{\pi}{2}p -(2k+1)b)}{(2k+1)^{1-p}}.  \nonumber
\end{equation}
\noindent
It turns out that this question is equivalent to validating the identity
\begin{equation}
(-1)^{k} 2^{-2n} \binom{2n+1}{n-k} (2k+1) = 
(2n+1) \sum_{j=k}^{n} (-1)^{j} 2^{-2j} \binom{n}{j} \binom{2j+1}{j-k} 
\label{sum99}
\end{equation}
\noindent
To this end, we employ the WZ-technology as explained in 
\cite{aequalsb}. This method produces the recurrence 
\begin{equation}
2(n+k+1)(n+1-k)u(n+1-k) - (n+1)(2n+3)u(n,k) = 0. 
\label{rec99}
\end{equation}
To prove (\ref{sum99}) simply check that both sides of (\ref{sum99}) 
satisfy the recurrence (\ref{rec99}) as well 
as the initial condition $u(0,0) = 1$. 

The identities 
\begin{eqnarray}
\ift x^{-p} \cos(x+b) \, dx & = & - \Gamma(1-p) \sin( b - \tfrac{p \pi}{2}) 
 \label{initcond} \\ 
\ift x^{-p} \sin(x+b) \, dx & = &  \Gamma(1-p) \cos( b - \tfrac{p \pi}{2}), 
\nonumber
\end{eqnarray}
\noindent
which are special cases of 
\begin{eqnarray}
\ift x^{-p} \cos(ax+b) \, dx  & = & -a^{p-1} \Gamma(1-p) 
\sin( b - \tfrac{p \pi}{2} ) \label{general} \\
\ift x^{-p} \sin(ax+b) \, dx  & = & a^{p-1} \Gamma(1-p) 
\cos( b - \tfrac{p \pi}{2} ), \nonumber
\end{eqnarray}
\noindent
show that the corresponding initial values in (\ref{gencos}) 
(respectively \ref{gensin}) match. The 
evaluations (\ref{general})  appear as 
$3.764.1$ and $3.764.2$ in \cite{gr}.  To establish (\ref{initcond}) 
expand $\cos(x+b)$ as $\cos x \cos b - \sin x \sin b$, use the 
change of variables $x \mapsto x^{p}$,  and Theorem 
\ref{main1}. 
\end{proof}

\medskip

We now discuss some definite integrals that follow from Theorem \ref{thmCS}. \\

\begin{example}
Differentiating (\ref{gencos}) with respect to $p$ and setting 
$p= \tfrac{1}{2}$ and $b=0$ gives, after the change of variables 
$x \mapsto x^{2}$, 
\begin{eqnarray}
& & \label{diff1} \\
\ift \log x \, \cos^{2n+1}x^{2} \, dx  & = & 
- \frac{\sqrt{\pi}}{2^{2n+3}} (\pi + 2 \gamma + 4 \log 2) 
\sum_{k=0}^{n} \binom{2n+1}{n-k} \frac{1}{\sqrt{4k+2}} \nonumber  \\
& - & \frac{\sqrt{\pi}}{2^{2n+2}} \sum_{k=0}^{n} \binom{2n+1}{n-k} 
\frac{\log(2k+1)}{\sqrt{4k+2}},
\nonumber
\end{eqnarray}
\noindent
where we have used the value $\Gamma'(1/2) = -\sqrt{\pi}( \gamma + 
2 \log 2)$. 
\end{example}

\begin{example}
Assume $ 0 < p, \, q < 1$. Multiplying (\ref{gencos}) by $b^{-q}$ and 
integrating over the half-line yields (after replacing $b$ by $y$)
\begin{eqnarray}
\ift \ift \frac{ \cos^{2n+1}(x+y) }{x^{p}y^{q}} dA  & =  & 
- \Gamma(1-p) \Gamma(1-q) \cos \left( \frac{\pi(p+q)}{2} \right) 
\nonumber \\
& \times & \sum_{k=0}^{n} \binom{2n+1}{n-k} \frac{(2k+1)^{p+q-2}}{2^{2n}}. 
\nonumber 
\end{eqnarray}
\noindent
In particular, for $n=0$, 
\begin{equation}
\ift \ift \frac{ \cos(x+y) }{x^{p}y^{q}} dA  =  
- \Gamma(1-p) \Gamma(1-q) \cos \left( \frac{\pi(p+q)}{2} \right). 
\end{equation}
\noindent
The derivative $\frac{\partial^{2}}{\partial p \partial q}$ at 
$p = q = \tfrac{1}{2}$ produces the evaluation
\begin{equation}
\ift \ift \frac{\log x \, \log y }{\sqrt{xy}} \cos(x+y) \, dx \, dy = 
( \gamma + 2 \log 2) \pi^{2} 
\end{equation}
\noindent
that we promised in the Introduction. 
\end{example}

\begin{example}
\label{last-example}
Iterating the method described in the previous example yields
\begin{equation}
\int_{{\mathbb{R}_{+}^{n}} } 
\left( \cos \| x \|^{2} \right) \cdot   \prod_{j=1}^{n} \log x_{j} \, dV = 
\frac{(-1)^{\Delta_{n}} \pi^{n/2}}{2^{2n}} \begin{cases}
           \realpart{\psi_{n}} \quad \text{ if } n \text{ is even}, \nonumber \\
           \imagpart{\psi_{n}} \quad \text{ if } n \text{ is odd}, \nonumber 
         \end{cases}
\end{equation}
\noindent
with 
\begin{equation}
\Delta_{n} = \frac{n(n+1)}{2}, \, \psi_{n} = \left( \gamma + 2 \log 2 + \frac{\pi i }{2} \right)^{n}
e^{\pi i n/4}.
\end{equation}
\noindent
Here $\| x \|^{2} = x_{1}^{2} + \cdots + x_{n}^{2}$ and $\gamma$ is Euler's constant. For instance, for 
$n=3$ we have
\begin{equation}
\ift \ift \ift \log x \log y \log z \cos( x^{2} + y^{2} + z^{2}) 
\, dx \, dy \, dz =
\frac{\pi^{3/2}}{8}( -16 \xi^{3} + 12 \xi^{2} \pi + 6 \xi \pi^{2} 
-  \pi^{3} ), \nonumber
\end{equation}
\noindent
where $\xi = \gamma + 2 \log 2$.
\end{example}

\bigskip

\noindent
{\bf Acknowledgements}. The third author acknowledges the 
partial support of NSF-DMS $040998$.  The second author was  partially
supported as a graduate student by the same grant.


\end{document}